\documentclass{article}
\usepackage{graphicx}
\usepackage{amsmath}
\usepackage{amsfonts}
\usepackage{amsthm}
\usepackage{amssymb}
\usepackage{color}
\usepackage{hyperref}
\usepackage{breakurl}

\begin{document}

\title{Odd One Out}

\author{Tanya Khovanova \\
\textit{Department of Mathematics, MIT}}

\maketitle


\section*{Introduction}

This article covers my second talk at the Gathering for Gardner in March,
2010. It is about an Odd One Out puzzle I invented, after having been inspired by Martin Gardner. I do not like Odd One Out questions; that is why I invented one.

\section*{A Standard Odd One Out Question}

Which one is the odd one out?

\begin{enumerate}
\item Cow
\item Hen
\item Pig
\item Sheep
\end{enumerate}

This is a standard question and the expected answer is: Hen, because it is a bird, while the others are mammals. Some people give ``wrong" answers, but when you ask them to explain, it gets interesting. For example, I've heard the following answers \cite{MyIQ}:

\begin{itemize}
\item Pig, because pig-meat is the only one that is non-kosher
\item Sheep, because it's the only word in the list that doesn't have three letters
\item Cow, because it's the only word in the list with letters in alphabetical order
\end{itemize}

Suppose you are going for a cup of coffee. Who would you rather invite along: the person who answers correctly or someone who comes up with a creative ``wrong" answer? Multiple-choice questions check standard reasoning, but punish non-triviality and originality.

\section*{Background}

Once I went out to dinner in Harvard Square with some puzzle lovers. One of
them had a Martin Gardner book and I scanned though it. On my way
home I realized that I had liked one of the puzzles in the book, but forgot the
details, so I reinvented it.

\section*{The Puzzle}

\begin{figure}[htp]
\begin{center}
\includegraphics[scale=0.8]{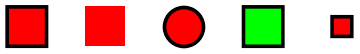}
\end{center}
\caption{Find the Odd One Out}
\label{odd}
\end{figure}

Which is the odd one out in Figure ~\ref{odd}?

\section*{Comments}

I published the puzzle on my blog \cite{OddOneOut}, and received a lot of comments. Some were just cute:

\begin{itemize}
\item The Odd One Out is number six, because it is invisible :-)
\item The Odd One Out is the observer: it's the only object that doesn't present any characteristic of the five objects shown.
\item None. Why should any of them be left out?
\item The right one because it is right.
\end{itemize}

Some were very interesting:

\begin{quote}
Intelligence is the fastest achievement of the simplest correct answer, not the ability to come up with the most lengthy and complicated one.
The simplest answer is: color difference is perceived before form or property set evaluation.
\end{quote}

Many people continued this theme. They gave more importance to one of the properties than to the others.

\begin{itemize}
\item Not Framed One --- it is the only one without a frame, and, in addition, it is the only one comprised of one color rather than two. So it differs by two features, as others differ only by one feature.
\item Circle One --- not only is it a different figure, but it also doesn't have any angles.
\end{itemize}

Will responded to a variety of the comments on the posting:

\begin{quote}
The prevailing opinion is that the first is the least unique and is therefore the oddest. But it is the mean and the others are one deviation from it. Can the mean be the statistical anomaly? Is that what the question asks for?
\end{quote}

To which Cedric replied:

\begin{quote}
Yes, I think the mean can be a statistical anomaly. The average person has roughly one testicle and one ovary, but a person with these characteristics would certainly be an anomaly.
\end{quote}

During my talk at the G4G9 someone suggested another reply:

\begin{quote}
A randomly chosen person almost surely will have more legs than the average number of legs per person.
\end{quote}

\section*{Meta}

Several comments were equating the meta property to other properties:

\begin{quote}
All right, the last four figures have special features; the first figure is special because it is normal. Hence, every figure is special and there are no odd ones here.
\end{quote}

It is a valid argument.
Many websites reprinted my puzzle. Some of them claimed that my intention was to create a puzzle where there is no odd one out. I have to disappoint them. That was not my intention. If I were designing the puzzle with no odd one out I would have designed the one in Figure ~\ref{noodd} (see \cite{NoOdd}):

\begin{figure}[htp]
\begin{center}
\includegraphics[scale=0.6]{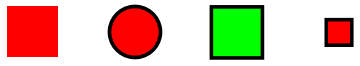}
\end{center}
\caption{Find the Odd One Out}
\label{noodd}
\end{figure}

During my talk someone pointed out that a simpler joke would be to
create the Odd One Out puzzle with several identical figures.

\section*{Martin Gardner}

\begin{figure}[htp]
\begin{center}
\includegraphics[scale=0.7]{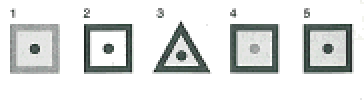}
\end{center}
\caption{Find the Odd One Out}
\label{Gardner}
\end{figure}

After the conference I decided to find Martin Gardner's book and take another look at the original problem. The book is titled \emph{The Colossal Book of Short Puzzles and Problems}. The puzzle is in Chapter 13, page 370, and was designed by Tom Ransom (See Figure ~\ref{Gardner}).

I was delighted to discover that Martin Gardner shares my dislike for Odd One Out questions. Here is what he wrote:

\begin{quote}
Sometimes there are so many different ways a symbol can be different from the others that brighter students are penalized by having to waste time deciding which symbol is ``most obviously different" to the person who designed the test.
\end{quote}

\section*{Another Puzzle}

In my G4G9 gift I included my original puzzle as well as the puzzle in Figure ~\ref{missing} that I designed with Gregory Bomash {see \cite{Missing}).

\begin{figure}[htp]
\begin{center}
\includegraphics[scale=0.8]{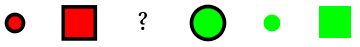}
\end{center}
\caption{What's Missing?}
\label{missing}
\end{figure}

\end{document}